\documentclass[12pt,a4paper]{amsart}
\usepackage{amssymb,eucal}

\newtheorem{theorem}{Theorem}
\newtheorem{conjecture}{Conjecture}
\newtheorem{lemma}{Lemma}
\newtheorem*{corollary}{Corollary}
\newtheorem*{definition}{Definition}
\newcommand{\set}[2]{\ensuremath{\{ #1 \>|\> #2 \}}}

\newcommand{\bigS}{\ensuremath{\text{\LARGE S}}}

\begin{document}

\title[Euler-Poincar\'e characteristic of a Lie superalgebra]
{How Euler would compute the Euler-Poincar\'e characteristic of a Lie superalgebra}
\author{Pasha Zusmanovich}
\address{
Department of Mathematics, Tallinn University of Technology, Ehitajate tee 5, 
19086 Tallinn, Estonia
}
\email{pasha.zusmanovich@ttu.ee}
\date{last revised August 8, 2011}
\thanks{\textsf{arXiv:0812.2255}; Expos. Math. \textbf{29} (2011), 345--360}

\begin{abstract}
The Euler-Poincar\'e characteristic of a finite-dimensional Lie algebra vanishes.
If we want to extend this result to Lie superalgebras, we should deal with infinite sums.
We observe that a suitable method of summation, which goes back to Euler,
allows to do that, to a certain degree.
The mathematics behind it is simple, we just glue the pieces of
elementary homological algebra, first-year calculus and pedestrian combinatorics 
together, and present them in a (hopefully) coherent manner.
\end{abstract}

\maketitle

\section*{Introduction}

Let $L$ be a finite-dimensional Lie algebra, $M$ an $L$-module. 
An \textit{Euler-Poincar\'e characteristic $\chi(L, M)$} is the alternating sum
\begin{equation}\label{sum-main}
\dim H^0(L, M) - \dim H^1(L,M) + \dim H^2(L,M) - ...
\end{equation}
As cohomology in dimension greater than $\dim L$ necessarily vanishes, the sum is 
finite.

The study of Euler-Poincar\'e characteristic of Lie algebras,
being probably less sophisticated than its differential-geometric, group-theoretic,
or category-theoretic counterparts, has nevertheless a relatively interesting history.

Apparently for the first time the question was considered in the Lie-algebraic 
framework by Chevalley and Eilenberg in \cite{ce}.
In that classical paper, cohomology rings of simple (classical) Lie algebras 
$\mathfrak g$
over algebraically closed field $K$ of characteristic zero with coefficients
in the trivial module $K$ were determined, and, in particular, it was noted that 
$\chi(\mathfrak g, K) = 0$ (op.~cit., p.~109). 
It was also proved that for any
finite-dimensional irreducible nontrivial $\mathfrak g$-module, the cohomology vanishes.
As every finite-dimensional $\mathfrak g$-module is completely reducible,
the cohomology with coefficients in any such module is a sum
of several copies of cohomology with coefficients in $K$, so the Euler-Poincar\'e 
characteristic for any finite-dimensional $\mathfrak g$-module vanishes too.

Chevalley and Eilenberg deduced their results from earlier results of 
Pontryagin and Hopf about Lie groups. Almost at the same time,
Koszul, also basing on Hopf's methods, gave a purely algebraic treatment 
of the subject in no less classical paper \cite{koszul}. 
In particular, he established that 
cohomology ring of a reductive Lie algebra over a field of characteristic zero
is isomorphic to the exterior algebra over the space of its primitive elements.
Although he did not mention the Euler-Poincar\'e characteristic 
explicitly, its vanishing immediately follows from the latter isomorphism.

Later, Goldberg \cite{goldberg} found an amazingly elementary proof of a more general 
result:
$\chi(L, K) = 0$ for any finite-dimensional Lie algebra $L$ over any field $K$.
(He considered only the case of the trivial module, but
the reasonings remain the same in the case of an arbitrary finite-dimensional module.
Also, it should be noted that he credited the idea of the proof in its full generality to 
the anonymous referee).
This anonymous referee's proof uses an elementary (and well-known) fact that 
the Euler-Poincar\'e characteristic of a finite cochain complex coincides with the 
Euler-Poincar\'e characteristic of its cohomology (referred sometimes as the 
Euler-Poincar\'e principle). 
We reproduce all ingredients of this proof
in an appropriate context later in the present note. 

Goldberg's paper apparently remained not widely known
(may be due to not very ``scientific'' place of publication), as later at least three 
papers appeared proving the vanishing of the Euler-Poincar\'e characteristic for different
classes of finite-dimensional algebras 
(that is, even not achieving the full generality): 
Malliavin-Brameret treated the cases of nilpotent \cite{malliavin-nilp} and then of 
solvable \cite{malliavin-solv} Lie algebras, and Pirashvili treated the cases of
Lie algebras over a field of characteristic zero and non-perfect algebras 
(over any field) \cite{pirashvili}. 
It is interesting that in all these papers, the same argument
as in Goldberg's paper was employed, but, so to say, ``in the wrong direction'': 
instead of passing from cohomology to the underlying cochain complex, they passed to the 
complex arising from the Hochschild-Serre spectral sequence, making arguments
more cumbersome and less general (this did not went unnoticed in the Zentralblatt review
0386.17006 of \cite{malliavin-solv} by Vogan).

What about Lie superalgebras?
An immediate difficulty arises: due to the presence of the odd part, the super analog
of the Chevalley-Eilenberg complex is infinite, so the sum (\ref{sum-main}) is,
in general, infinite and divergent. 
Even if we restrict considerations to Lie superalgebras and modules whose cohomology
vanishes for a large enough degree (this is a far from trivial and fairly
interesting class of Lie superalgebras, as was demonstrated in \cite{gruson}),
the underlying cochain complex is still infinite, and we can not pass to it easily
as in the ordinary Lie algebras case.
As Borcherds succinctly put it in \cite{borcherds}: ``Strictly speaking... 
the alternating sum (\ref{sum-main}) does not make sense unless the groups $\bigwedge^n$ 
are all finite dimensional and almost all $0$. I will deal with this problem by ignoring 
it''.

Of course, people successfully applied the Euler-Poincar\'e principle to 
infinite (co)chain complexes. However, in all such cases, that or another sort of 
finiteness condition was present, which allowed to reduce the situation to the usual case 
of finite complexes.
For example, in, perhaps, the most celebrated such application --
proof of the Weyl-Kac character formula and various combinatorial identities by 
homological methods due to Garland and Lepowsky (see \cite[Chapter 3, \S\S 2-3]{fuchs})
-- the underlying Chevalley-Eilenberg complex computing the homology of the ``nilpotent''
part of a Kac-Moody algebra has a structure of a semisimple module over a Cartan subalgebra
with finite-dimensional weight spaces, what allows to decompose the whole complex into
the (infinite) direct sum of finite ones.

In absence of any finiteness conditions, we apparently do not have another choice
than to employ an appropriate summation method
of infinite series. With a suitably modified definition
of the Euler-Poincar\'e characteristic we are able to carry on the Goldberg's 
proof in the super case, under some assumption of convergence
and dimensions of the homogeneous components.

In fact, we consider a more general setting of color Lie algebras.
Though, compared to the super case, we are able to give even a less complete answer,
it has an interesting conceptual feature that the notions of 
dimension and  the Euler-Poincar\'e characteristic, are, essentially, the same.

Our results could be interpreted in the following interesting way: if 
the zero component of a color Lie algebra is large enough (i.e., if a color algebra is 
``close'' to an ordinary Lie algebra), its Euler-Poincar\'e characteristic
behaves like in the ordinary Lie algebra case.
On the contrary, and this is peculiar to the super/color situation, if dimensions of the zero and odd components coincide, 
the Euler-Poincar\'e characteristic does not vanish!

The next two short sections contain recapitulations on all necessary definitions
and results from the color algebra, \S 3 contains definition of the Euler-Poincar\'e
characteristic suitable for the color case based on the Abel summation, 
as well as the needed preparatory statements about it, 
and \S 4 contains formulation and proof of the main theorem.
In \S 5 we speculate about the possibility to use other summation methods.
As the goal of this note is to say as much as possible about the Euler-Poincar\'e
characteristic, entirely avoiding at the same time to look at cohomology itself, 
there is no shortage of further questions which may require more sophisticated approaches. 
In the last \S 6 we point out some of them.

We choose to argue about cohomology, but all the definitions, reasonings and results
could be carried over to homology in an obvious way.

\section{Recapitulation on color linear algebra}\label{section-dim}

Here we record some elementary facts about color dimension
in the generality suited for our needs. 
Essentially, they could be found in different variations and using different terminology
in many places in the literature (see, for example, \cite{kang} and references therein). 
Otherwise, we follow \cite{scheunert} and \cite{scheunert-zhang}, where all other 
necessary details may be found.

Let $G$ be an abelian group (written additively) and $K$ be a field. 
We work in the category
of $G$-color vector spaces over $K$, i.e., $G$-graded vector spaces 
with suitably defined ($G$-graded)
morphisms, subspaces, tensor products, etc. 
Most of the time we are speaking, by abuse of language, on just 
\textit{color vector spaces}, assuming it is clear from the context
which group $G$ is meant.

To circumvent the standard set-theoretic issues (like nonexistence of the set of all
cardinal numbers), assume that all elements and sets under consideration are lying in some 
(infinite) universe,
and let $\mathcal C$ be a ring consisting of all cardinal numbers 
less than the cardinality of the universe,
as well as their formal negations, with obvious operations of addition and multiplication.
(This, however, is not of paramount importance for our further considerations, as we are
interested primarily in the finite-dimensional case, in which the role of $\mathcal C$
is played by $\mathbb Z$).

\begin{definition}
A ($G$-color) dimension of a $G$-color vector space 
$V = \bigoplus_{\alpha\in G} V_\alpha$, denoted as $\dim_G V$, is an element of the group 
ring $\mathcal C[G]$ equal to
\begin{equation}\label{formula-sum}
\sum_{\alpha\in G} (\dim V_\alpha) e^\alpha,
\end{equation}
where $\dim V_\alpha$ is an ordinary dimension of vector spaces over $K$.
$V$ is said to be finite-dimensional if all graded components 
$V_\alpha$ are finite-dimensional (irrespective of whether $G$ is finite or 
not)\footnote[2]{
In some specific situations, for example, when the vector space $V$ is a highest weight
module over a simple classical or Kac-Moody Lie algebra, with grading induced by 
the action of the Cartan subalgebra, and the group $G$ is generated
by the corresponding roots, what we call a color dimension here, is called a \textit{character} 
(cf. the famous Weyl-Kac character formula).
}. 
\end{definition}

(Note that $e^\alpha$'s in formula (\ref{formula-sum}) are merely formal symbols 
facilitating passage from the additive notation in $G$ to the multiplicative notation in 
$\mathcal C[G]$, i.e., satisfying the usual properties of 
the exponent: $e^0 = 1 \in \mathcal C$, $e^\alpha e^\beta = e^{\alpha + \beta}$).

In what follows, we will always specify subscript $G$ in $\dim_G$, to distinguish it from 
the ordinary dimension $\dim$ (which is, of course, is just the color dimension in the case
where $G$ is the trivial group).

Obviously, dimension of a finite-dimensional $G$-color vector space
lies in the subring $\mathbb Z[G]$ of $\mathcal C[G]$.

The color dimension enjoys the usual additivity and multiplicativity 
properties of the ordinary dimension:

\begin{lemma}\label{lemma-dimension}\hfill
\begin{enumerate}
\item
If $0 \to V \to U \to W \to 0$ is a short exact sequence of color vector spaces,
then $\dim_G U = \dim_G V + \dim_G W$. In particular (or rather, equivalently), 
$\dim_G (V \oplus W) = \dim_G V + \dim_G W$.
\item $\dim_G (V\otimes W) = \dim_G V \cdot \dim_G W$.
\item If $V$ is finite-dimensional, then $\dim_G Hom(V, W) = \dim_G V \cdot \dim_G W$.
\end{enumerate}
\end{lemma}

\begin{proof}
This is obvious and reflects addition and multiplication rules in the group ring.
\end{proof}

The following can be considered as a generalization of the formula (\ref{formula-sum}):

\begin{lemma}\label{lemma-direct-product}
Let $G$, $H$ be two abelian groups, 
$V = \bigoplus_{(\alpha, \beta)\in G\times H} V_{\alpha, \beta}$
be a $G \times H$-color vector space,
$\Phi: \mathcal C [G \times H] \simeq \mathcal C[G][H]$ be the canonical isomorphism.
Then
$$
\Phi(\dim_{G\times H} V) = \sum_{\beta \in H} (\dim_G V_\beta) e^\beta ,
$$
where
$V_\beta$ is a $G$-color vector space equal to 
$\bigoplus_{\alpha\in G} V_{\alpha, \beta}$.
\end{lemma}

\begin{proof}
Amounts to writing the double sum 
$\bigoplus_{\alpha\in G, \beta\in H} V_{\alpha, \beta}$ in two ways, by grouping
elements by $G$ and by $H$.
\end{proof}

In a more general situation, when there is a surjective group homomorphism 
$\varphi: G \to H$, a $G$-color vector space $V$ can be considered as an 
$H$-color vector space, by grouping together elements lying in the 
same conjugacy class: $V = \bigoplus_{\beta \in H} V_\beta$,
$V_\beta = \bigoplus_{\alpha \in \varphi^{-1}(\beta)} V_\alpha$.
Accordingly, we may consider the dimension $\dim_H V$ of $V$ as an $H$-color vector space,
which coincides with the image of $\dim_G V$ under the induced homomorphism of group
rings $\mathcal C[G] \to \mathcal C[H]$.

In particular, each color vector space has an \textit{ordinary dimension} as an
ordinary vector space, which is just the sum of dimensions of all its homogeneous 
components.

\section{Recapitulation on color Lie algebras and their cohomology}\label{section-color}

Though the number of publications on color algebras estimated nowadays
by hundreds (and on superalgebras by thousands), there is still unfortunate ambiguity in 
terminology and notations, so we briefly set them up in this section. 
We continue to follow \cite{scheunert} and \cite{scheunert-zhang}.

Let $G$ be, like in \S \ref{section-dim}, an abelian group. 
A \textit{commutation factor} on $G$ is a map $\varepsilon: G \times G \to K$ such that
\begin{gather}
\varepsilon(\alpha, \beta) \varepsilon(\beta, \alpha) = 1 \label{one} \\
\varepsilon(\alpha + \beta, \gamma) = 
\varepsilon(\alpha, \gamma) \varepsilon(\beta, \gamma)    \label{two} \\
\varepsilon(\alpha, \beta + \gamma) = 
\varepsilon(\alpha, \beta) \varepsilon(\alpha, \gamma)    \label{three}
\end{gather}
for all $\alpha, \beta, \gamma \in G$.
A \textit{$(G, \varepsilon)$-color Lie algebra} is a color vector space 
$L = \bigoplus_{\alpha\in G} L_\alpha$
equipped with multiplication $[\>\cdot\>, \>\cdot\>]$ which preserves $G$-grading and
satisfies $\varepsilon$-skew-symmetricity and $\varepsilon$-Jacobi identity
(in fact, conditions (\ref{one})-(\ref{three}) defining a commutation factor
may be deduced from the preservation of $G$-grading under $[\>\cdot\>, \>\cdot\>]$ 
and $\varepsilon$-skew-symmetricity).

Again, by abuse of language, we are talking about just \textit{color Lie algebras},
assuming the group $G$ and the commutation factor $\varepsilon$
are clear from the context. The important case $G = \mathbb Z_2$ and 
\begin{equation}\label{comm-factor}
\varepsilon(\alpha, \beta) = (-1)^{\alpha\beta}
\end{equation}
warrants a special name of a \textit{Lie superalgebra}.
The underlying vector space of a Lie superalgebra (i.e., a $\mathbb Z_2$-color vector space)
is called \textit{superspace}.
The dimension of a Lie superalgebra (and, more general, of a superspace) 
is an element of the ring of dual numbers $\mathbb Z[\mathbb Z_2] \simeq
\mathbb Z[\theta]/(\theta^2 - 1)$, and often is denoted as a pair $(n | m)$, where
$n$ and $m$ are dimensions of the even and odd parts, respectively.

It follows from (\ref{one}) that $\varepsilon(\alpha, \alpha) = \pm 1$ for all
$\alpha\in G$. This means that there is a homomorphism $G \to \mathbb Z_2$
with the kernel $G_{even} = \set{\alpha\in G}{\varepsilon(\alpha, \alpha) = 1}$
and non-trivial conjugacy class 
$G_{odd} = \set{\alpha\in G}{\varepsilon(\alpha, \alpha) = -1}$, and
$G = G_{even} \cup G_{odd}$. Note that the order $|\alpha|$ of each element 
$\alpha\in G_{odd}$ is even, provided it is finite.
Denoting $L_i = \bigoplus_{\alpha\in G_i} L_\alpha$, $i \in \{even, odd\}$, we have
$L = L_{even} \oplus L_{odd}$.  
Subspaces $L_{even}$ and $L_{odd}$ (in fact, $L_{even}$ is a subalgebra) 
are called \textit{even part} and \textit{odd part}
of $L$ respectively\footnote[2]{
In many texts the even and odd parts of a color Lie algebra $L$
are denoted as $L_0$ and $L_1$ respectively. In our notations $L_0$ is a genuine 
``zero part'', a homogeneous component corresponding to the zero element $0\in G$. 
Of course, when $L$ is a superalgebra, $L_0 = L_{even}$, but in general $L_0$ is a 
direct summand of $L_{even}$.}. 
Thus, each color Lie algebra can be considered as a superspace
(but, in general, not as a Lie superalgebra), and has a corresponding dimension, 
which will be called \textit{superdimension} and denoted as $\dim_{\mathbb Z_2}$. 
Clearly, by what is written at the end of
the previous section, the superdimension is obtained from the color dimension by 
summing up dimensions of all odd and even homogeneous graded components respectively.

There is straightforward notion of a color module $M$ over a color Lie algebra $L$
(the action of $L$ on $M$ is compatible with $G$-grading).
A \textit{cohomology} of a color Lie algebra $L$ 
with coefficients in the color module $M$
is defined as cohomology of the \textit{color Chevalley-Eilenberg complex}
$$
C^0(L,M) \to C^1(L,M) \to C^2(L,M) \to \dots
$$
where $C^n(L,M)$ is the $G \times \mathbb N$-color space of $n$-fold multilinear 
$\varepsilon$-skew-symmetric maps from $L$ to $M$ (the $G$-grading is inherited from $L$ and $M$,
and the $\mathbb N$-grading is the usual cohomology degree), 
and differential is defined by the long formula
containing a lot of confusing signs, generalizing the differential in the ordinary
Chevalley-Eilenberg cochain complex (whose exact appearance is inessential to us here
due to wonders of invariance of the Euler-Poincar\'e characteristic under the
operation of taking cohomology, see Lemma \ref{principle} below). 
The space $C^n(L,M)$ may be canonically identified
with $Hom(\bigwedge_\varepsilon^n L, M)$, where $\bigwedge_\varepsilon^n L$ is the
$n$-fold $\varepsilon$-exterior power of the color space $L$.
Note that the latter is also a $G$-color vector space, 
the grading is defined by adding the gradings of all factors.

One has
\begin{equation}\label{epsilon-isom}
\bigwedge\nolimits_\varepsilon^*(V \oplus W) \simeq 
(\bigwedge\nolimits_\varepsilon^* V) \>\widehat\otimes\> 
(\bigwedge\nolimits_\varepsilon^* W)
\end{equation}
as $G \times \mathbb N$-color vector spaces. 
Here $\widehat\otimes$ denotes the tensor product twisted by a commutation factor
which is obtained by combining the commutation factor $\varepsilon$ with another one
obtained by the standard formula (\ref{comm-factor}) from the $\mathbb N$-grading.
Note, however, that this twisting obviously does not affect the dimension of color 
spaces to be tensored, i.e.
\begin{equation}\label{twist}
\dim_{G \times \mathbb N} (V \otimes W) = \dim_{G \times \mathbb N} (V \>\widehat\otimes\> W)
\end{equation}
for any two color vector spaces $V$ and $W$.
Applying the isomorphism (\ref{epsilon-isom}) to the whole color Lie algebra
$L = \bigoplus_{\alpha\in G} L_\alpha$, and noting that the $\varepsilon$-exterior power
coincides with the ordinary exterior power and with the ordinary symmetric power
for even and odd components respectively, one get the isomorphism of 
$G \times \mathbb N$-color vector spaces:
\begin{equation}\label{decomp}
\bigwedge\nolimits_\varepsilon^* L \simeq 
\widehat\bigotimes_{\alpha\in G_{even}} \bigwedge\nolimits^*(L_\alpha) \>\widehat\otimes\>
\widehat\bigotimes_{\alpha\in G_{odd}} \bigS^*(L_\alpha).
\end{equation}

\section{The Euler-Poincar\'e characteristic}\label{section-char}

Let $V = \bigoplus_{n\ge 0} V_n$ be a (possibly infinite-dimensional) 
$\mathbb N$-graded $G$-color vector space (i.e. each homogeneous component
$V_n = \bigoplus_{\alpha\in G} V_{\alpha,n}$ is a $G$-color vector space), 
all whose graded components are finite-dimensional. 
It can be considered as a $G \times \mathbb N$-color vector space 
$V = \bigoplus_{(\alpha, n) \in G\times \mathbb N} V_{\alpha, n}$.

Consider $\dim_{G\times \mathbb N} V$. Applying Lemma \ref{lemma-direct-product}
and writing $\mathbb N$ as multiplicative semigroup generated by a single element
$t$ of infinite order, one can write this dimension as a formal power series 
\begin{equation}\label{formal}
\sum_{n\ge 0} (\dim_G V_n)\, t^n.
\end{equation}
It lies in $\mathbb Z[G][t]$,
denoted as $\chi_G (V,t)$ and is called a \textit{Poincar\'e series}\footnote[2]{
Sometimes in the literature one calls this \textit{Hilbert series},
leaving the term \textit{Poincar\'e series} for the specific situation of
cohomology space graded by cohomology degree.} 
of $V$.

Now we want to define the Euler-Poincar\'e characteristic of a color space $V$, having
the cohomology of color Lie algebras in mind. 
A (naive) Lie-algebraic approach would be to define it as a series 
\begin{equation}\label{sum}
\sum_{n\ge 0} (-1)^n \dim_G V_n ,
\end{equation}
obtained by substituting 
$t = -1$ in the Poincar\'e series. When $V$ is the cohomology of a finite-dimensional
Lie algebra $L$, this is all right: due to anticommutativity of $L$ the standard
Chevalley-Eilenberg complex is finite, and the sum is finite. 

This is not so, however,
for Lie superalgebras, due to presence of the odd part: the cochain complex is infinite,
so one could try to make sense of (\ref{sum}) by considering its convergence in
the group ring $\mathbb R[G]$.
We will understand the convergence of sequences (and series) 
with elements in $\mathbb R[G]$ (considered for that
matter as the set of functions with a finite support on $G$)
in the sense of the weak convergence, what is equivalent to the 
$G$-componentwise convergence.
Since all terms of the series lie in $\mathbb Z[G]$, the convergence
implies that for each $\alpha\in G$,
$V_{\alpha,n}$ vanishes for all but a finite number of $n$'s.
This does not necessarily mean that 
$V_{\alpha,n}$ vanish uniformly over $G$, but this is of course so if $G$ is finite,
so in that case the convergence of (\ref{sum}) implies the finite-dimensionality of $V$.
This is anything but interesting as it breaks our scheme below of passing from cohomology
to the underlying cochain complex.

The situation would not be much eased even if we consider in (\ref{sum}) the ordinary 
dimension instead of the color dimension: the series 
\begin{equation}\label{sum-r}
\sum_{n\ge 0} (-1)^n \dim V_n ,
\end{equation}
still diverges (in $\mathbb R$) unless $V_n$ is zero for all but a finite number of $n$'s.

Thus, the series (\ref{sum}) or (\ref{sum-r}) diverge at least in the most interesting cases, and
to make a sense of them, one should apply one of the summation methods of divergent 
series.

We choose the Abel summation, for three reasons. First, to justify the title of this note.
Consider the simplest example of the $(0|1)$-dimensional Lie superalgebra 
(zero even component, $1$-dimensional odd component, and the trivial multiplication)
with the trivial $(1|0)$-dimensional module. The coboundaries are zero, so cohomology
coincides with the whole cochain complex, which has dimension $(1|0)$ in each even
component and dimension $(0|1)$ in each odd component.
The infinite sum (\ref{sum-r}) in this situation becomes 
$$
1 - 1 + 1 - 1 + \dots
$$
But it was Euler who first successfully and methodically computed sums like this, and it is 
merely a matter of chance that this particular summation method is named after Abel
(who considered the idea of summation of divergent series with disgust)\footnote[2]{
Of course, Euler got his share of fame also here, and there is another summation method named 
after him.
}.

Second, the Abel summation method is powerful enough, at least it is more powerful
than the most popular and simple methods.

Third, this method makes a perfect connection with the formal series (\ref{formal}),
so the Euler-Poincar\'e characteristic, like in the finite Lie-algebraic situation, 
still obtained in a very straightforward way from the Poincar\'e series, 
which, in its turn, is nothing but a properly understood color dimension\footnote[3]{
There is another popular summation method, more powerful than the Abel summation,
and still having advantage to be intimately related to the formal series (\ref{formal}) -- 
by analytic continuation:
if we consider  the parameter $t$ in (\ref{formal}) as lying in $\mathbb C$, 
and the radius of convergence of this series happens to be non-zero (but possibly $<1$), 
then it defines an analytic function in some 
circle which, perhaps, could be analytically prolonged to the point $t=-1$. 
However, all the (infinite) series we are dealing with here have the radius of
convergence equal to $1$, and in such case these two methods are equivalent.
See question 2 in \S \ref{section-questions}.
}.

So, following Euler (and, more lately, Berger, Leinster, Propp and others,
who applied summation of divergent Euler-Poincar\'e-like 
series to other situations arising in geometry and category theory -- 
see, for example, \cite{berger-l} and references therein), we will cheat a bit by 
considering a formal parameter $t$ in (\ref{formal}) as lying in $\mathbb R$, 
and considering 
\begin{equation}\label{limit}
\lim_{t \to -1+0} \sum_{n\ge 0} (\dim_G V_n)\, t^n.
\end{equation}
If exists, this limit, lying in $\mathbb R[G]$, will be called
the \textit{Euler-Poincar\'e characteristic} of a color space $V$ and will be denoted
as $\chi_G(V)$.

If, in (\ref{formal}) and (\ref{limit}), instead of $G$-dimension we will consider
ordinary or super dimension, we will arrive at the notion of \textit{ordinary} or 
\textit{super Poincar\'e series} and \textit{Euler-Poincar\'e characteristic} of a color 
vector space $V$, denoted as $\chi(V,t)$, $\chi_{\mathbb Z_2}(V,t)$ and 
$\chi(V)$, $\chi_{\mathbb Z_2}(V)$, respectively
(as in the case of $\dim$, we omit the subscript in the ordinary case, when the group
$G$ is trivial).
Note that the existence of $\chi_G(V)$ implies the existence of $\chi_{\mathbb Z_2}(V)$ 
which, in turn, implies the existence of $\chi(V)$, but the opposite implications are, 
in general, not true: for example, for the just considered cochain complex of 
the abelian $(0|1)$-dimensional Lie superalgebra with coefficients in the trivial 
$(1|0)$-dimensional module, the ordinary Euler-Poincar\'e characteristic is equal to 
$\frac{1}{2}$, while the super Euler-Poincar\'e characteristic does not exist: 
the super Poincar\'e series has the form $\frac{1 - t\theta}{1 - t^2}$.

A \textit{Poincar\'e series} and \textit{Euler-Poincar\'e characteristic} of a color Lie 
algebra $L$ and an $L$-module $M$, denoted as $\chi_G(L,M,t)$ and $\chi_G(L,M)$ respectively, 
is defined as Poincar\'e series and Euler-Poincar\'e characteristic
of a cohomology space $H^*(L,M)$ (considered with the standard grading by cohomology
degree). As noted above, in the case of an ordinary Lie algebra $L$, the series involved
degenerate to finite sums, and our definition coincides with the standard one.
The ordinary and super variants $\chi(L,M,t)$, $\chi_{\mathbb Z_2}(L,M,t)$ and 
$\chi(L,M)$, $\chi_{\mathbb Z_2}(L,M)$ are defined analogously.

\begin{lemma}\label{lemma-one-plus-t}
Let $C^*$ be a color cochain complex all whose graded components 
are finite-dimensional, and $B^*$ be the space of its coboundaries.
Then 
\begin{equation}\label{eq-one-plus-t}
\chi_G(C^*,t) - \chi_G(H(C^*),t) = (1+t) \chi_G(B^*,t).
\end{equation}
\end{lemma}

\begin{proof}
Denote, additionally, by $Z^*$ the space of cocycles.
There are two short exact sequences of $G \times \mathbb N$-color vector spaces:
\begin{gather}
0 \to Z^* \to C^* \to B^*[1] \to 0 \label{exact-1} \\
0 \to B^* \to Z^* \to H(C^*) \to 0 \label{exact-2}
\end{gather}
($B^*[1]$ is a graded space obtained from $B^*$ by shifting grading to 1 forwards).

By Lemma \ref{lemma-dimension}(i), 
\begin{align*}
\chi_G(C^*,t) &= \chi_G(Z^*,t) + t\chi_G(B^*,t) \\
\chi_G(Z^*,t) &= \chi_G(B^*,t) + \chi_G(H(C^*),t)
\end{align*}
what yields the desired equality.

This essentially constitutes the first part of the proof in \cite{goldberg},
and is a well-known argument in homological algebra (see, for example, \cite[\S 2.8]{bourbaki}).
\end{proof}

\begin{lemma}\label{principle}
Let $C^*$ be a color cochain complex all whose graded components are finite-dimensional,
and $B^*$ be the space of its coboundaries.
If both $\chi_G(C^*)$ and $\chi_G(B^*)$ exist,
then $\chi_G(H(C^*))$ exists and coincides with $\chi_G(C^*)$.
\end{lemma}

\begin{proof}
This amounts to ``substituting'' $t=-1$ into (\ref{eq-one-plus-t}). However, as 
we are dealing with infinite series, 
the actual arguments will be a bit more careful as follows.

The existence of $\chi_G(C^*)$ implies that the radius of convergence of the series
$\sum_{n\ge 0} (\dim C_{\alpha}^n)t^n$ for each $\alpha \in G$ is $\ge 1$. Thus
$$
1 \le 
\frac{1}{\limsup_{n \to \infty} \sqrt[n]{\dim C_{\alpha}^n}} \le 
\frac{1}{\limsup_{n \to \infty} \sqrt[n]{\dim H_{\alpha}^n}}
$$
which shows that the radius of convergence of the series 
$\sum_{n\ge 0} (\dim H_{\alpha}^n)t^n$ is also $\ge 1$.
Hence, for each $|t| < 1$, the series $\sum_{n\ge 0} (\dim_G H^n)t^n$ 
converges to some element from $\mathbb R[G]$.
Now assuming in (\ref{eq-one-plus-t}) $t\in \mathbb R$ and passing on both sides
to $\lim_{t\to -1+0}$, we get
$$
\chi_G(C^*) - \chi_G(H(C^*)) = \lim_{t \to -1+0}(1+t) \> \chi_G(B^*) = 0.
$$
The existence of $\chi_G(C^*)$ and $\chi_G(B^*)$ implies the existence of $\chi_G(H(C^*))$,
and the statement of the lemma follows.
\end{proof}

In what follows, $T^*(V)$ denotes, as usual, the tensor algebra generated by a vector space $V$,
with a standard $\mathbb N$-grading.

\begin{lemma}\label{multisec}
Let $G$ be a finite abelian group, $V$ be a $G$-color vector space whose dimension is 
concentrated in one degree $e^\alpha$, and $W$ be a subspace 
(as $G \times \mathbb N$-color vector space) of $T^*(V)$. Then
\begin{equation*}
\chi_G(W,t) = \frac{1}{|\alpha|} 
\sum_{i=0}^{|\alpha| - 1} \Big( 
\sum_{j=0}^{|\alpha| -1} \frac{\chi(W, \omega^j t)}{\omega^{ij}} \Big) e^{i\alpha},
\end{equation*}
where $\omega$ is a primitive $|\alpha|$th root of unity.
\end{lemma}

\begin{proof}
Let $W = \bigoplus_{n\ge 0} W_n$ be decomposition of $W$ as a subspace of $T^*(V)$ 
with respect to $\mathbb N$-grading.
Since dimension of the $n$-fold tensor product $V \otimes \dots \otimes V$ is concentrated
in one degree $e^{n\alpha}$, this decomposition yields
$$
\chi_G(W,t) = \sum_{n=0}^\infty (\dim W_n) \> e^{n\alpha} t^n
$$
and
$$
\chi(W, t) = \sum_{n=0}^\infty (\dim W_n) \> t^n .
$$
Since $|\alpha|$ is finite and different from zero, $\chi_G(W, t)$ could be rewritten as
$$
\sum_{i=0}^{|\alpha| - 1} \Big( \sum_{n=0}^\infty (\dim W_{kn+i}) \> t^{kn+i} \Big) 
e^{i\alpha} ,
$$
the coefficients of $e^{i\alpha}$ being exactly multisections (\cite[\S 4.3]{riordan})
of the power series $\chi(W,t)$, what implies the desired formula.
\end{proof}

To formulate some of the next technical results, we will need the following auxiliary
notation.
Fix an element $\alpha\in G$ of a finite order and consider a 
finite-dimensional $G$-color vector space $W$ having dimension $1$ in each degree 
$e^{i\alpha}$. The Euler-Poincar\'e characteristic of $W$ is equal to 
$$
\sum_{i=0}^{|\alpha| - 1} (-1)^i e^{i\alpha} ,
$$
and its ``scaled'' version
$$
\frac{1}{|\alpha|} \sum_{i=0}^{|\alpha| - 1} (-1)^i e^{i\alpha}
$$
will be denoted as $\overline\chi(\alpha)$.
Like the ``usual'' Euler-Poincar\'e characteristic, it lies in the group ring
$\mathbb Z[G]$. 

\begin{lemma}\label{complex}
Let $L$ be a finite-dimensional color Lie algebra.
\begin{enumerate}
\item
$\chi (\bigwedge\nolimits_\varepsilon^* L)$ exists and is equal to
$$
\begin{cases}
0,                          & L_{even} \ne 0  \\
\frac{1}{2^{\dim L_{odd}}}, & L_{even} = 0.
\end{cases}
$$

\item
$\chi_{\mathbb Z_2} (\bigwedge\nolimits_\varepsilon^* L)$ exists if and only if 
$\dim L_{even} \ge \dim L_{odd}$ and is equal to
$$
\begin{cases}
0,                                            & \dim L_{even} > \dim L_{odd}  \\
\big( \frac{1 - \theta}{2} \big)^{|G_{odd}|}, & \dim L_{even} = \dim L_{odd}.
\end{cases}
$$

\item
$\chi_G (\bigwedge\nolimits_\varepsilon^* L)$ exists if and only if 
$\dim L_0 \ge \dim L_{odd}$ and is equal to 
$$
\begin{cases}
0,                                               & \dim L_0 > \dim L_{odd}  \\
\prod_{0\ne \alpha\in G_{even}} (1 - e^\alpha)^{\dim L_\alpha}
\prod_{\alpha\in G_{odd}} \overline\chi(\alpha), & \dim L_0 = \dim L_{odd}.
\end{cases}
$$
\end{enumerate}
\end{lemma}

\begin{proof}
In view of (\ref{decomp}), (\ref{twist}) and Lemma \ref{lemma-dimension}(ii),
\begin{equation}\label{prod}
\chi_G(\bigwedge\nolimits_\varepsilon^* L, t) =
\prod_{\alpha\in G_{even}} \chi_G(\bigwedge\nolimits^*(L_\alpha), t) 
\prod_{\alpha\in G_{odd}}  \chi_G(\bigS^*(L_\alpha), t),
\end{equation}
same for $\chi$ and $\chi_{\mathbb Z_2}$.
Each complex $\bigwedge\nolimits^*(L_\alpha)$ is finite, hence its Poincar\'e
series reduces to the finite sum 
\begin{equation}\label{skew-finite}
\chi_G(\bigwedge\nolimits^*(L_\alpha), t) = 
\sum_{n=0}^{\dim L_\alpha} \dim \bigwedge\nolimits^n(L_\alpha) \> e^{n\alpha} t^n =
\sum_{n=0}^{\dim L_\alpha} \binom{\dim L_\alpha}{n} \> e^{n\alpha} t^n = 
(1 + te^\alpha)^{\dim L_\alpha},
\end{equation}
and, since $\alpha\in G_{even}$,
\begin{equation}\label{ord-skew}
\chi(\bigwedge\nolimits^*(L_\alpha), t) = 
\chi_{\mathbb Z_2} (\bigwedge\nolimits^*(L_\alpha), t) = 
(1+t)^{\dim L_{\alpha}}.
\end{equation}
However, for the reason which should be clear in a few seconds, we will treat this 
expression also in a bit more cumbersome way.
Since $L$ is finite-dimensional, $G$ is finite, so applying Lemma \ref{multisec} we get
\begin{equation}\label{color-skew}
\chi_G(\bigwedge\nolimits^*(L_\alpha), t) = \frac{1}{|\alpha|} 
\sum_{i=0}^{|\alpha| - 1} \Big( 
\sum_{j=0}^{|\alpha| - 1} \frac{(1 + \omega_\alpha^j t)^{\dim L_\alpha}}{\omega_\alpha^{ij}}
\Big) e^{i\alpha} , 
\end{equation}
where $\omega_\alpha$ is a primitive $|\alpha|$th root of unity.

Similarly, 
\begin{equation}\label{symm}
\chi_G(\bigS^*(L_\alpha), t) = 
\sum_{n=0}^\infty \dim \bigS^n (L_\alpha) \> e^{n\alpha} t^n =
\sum_{n=0}^\infty \binom{\dim L_\alpha + n-1}{n} \> e^{n\alpha} t^n ,
\end{equation}
what implies 
\begin{equation}\label{ord-symm}
\chi(\bigS^*(L_\alpha), t) = \sum_{n=0}^\infty \binom{\dim L_\alpha + n-1}{n} t^n =
\frac{1}{(1 - t)^{\dim L_\alpha}} ,
\end{equation}
and application of Lemma \ref{multisec} yields
\begin{equation}\label{color-symm}
\chi_G(\bigS^*(L_\alpha), t) = 
\frac{1}{|\alpha|} \sum_{i=0}^{|\alpha| - 1} \Big( 
\sum_{j=0}^{|\alpha| - 1} \frac{1}{\omega_\alpha^{ij} (1 - \omega_\alpha^j t)^{\dim L_\alpha}}
\Big) e^{i\alpha} 
\end{equation}
and
\begin{multline}\label{super-symm}
\chi_{\mathbb Z_2} (\bigS^*(L_\alpha), t) \\ = \frac{1}{2} 
\Big( 
\frac{1}{(1 - t)^{\dim L_\alpha}} +
\frac{1}{(1 + t)^{\dim L_\alpha}} 
\Big) 
+ \frac{1}{2}
\Big( 
\frac{1}{(1 - t)^{\dim L_\alpha}} -
\frac{1}{(1 + t)^{\dim L_\alpha}} 
\Big) \theta .
\end{multline}
Now, putting (\ref{prod}), (\ref{ord-skew}) and (\ref{ord-symm}) together, we get:
\begin{equation}\label{chi-ord}
\chi (\bigwedge\nolimits_\varepsilon^* L, t) = 
\frac
{(1+t)^{\sum_{\alpha\in G_{even}} \dim L_\alpha}}
{(1-t)^{\sum_{\alpha\in G_{odd}} \dim L_\alpha}} =
\frac
{(1+t)^{\dim L_{even}}}
{(1-t)^{\dim L_{odd}}} ,
\end{equation}
and part (i) follows.
This is, essentially, the second part of the arguments in \cite{goldberg}
(somewhat obscured by the presence of odd components), and finishes
the ordinary Lie algebra case.
The formula (\ref{chi-ord}) was also noted in \cite[Proposition 12]{piont-s}.

Similarly, collecting (\ref{prod}), (\ref{ord-skew}) and (\ref{super-symm}), we get a 
formula for $\chi_{\mathbb Z_2} (\bigwedge\nolimits_\varepsilon^* L, t)$, which is, 
however, more cumbersome.
Fortunately, we do not need to care about the exact expression, 
but merely note that it could be written as the sum of terms of the form 
$$
\pm \frac{1}{2^{|G_{odd}|}} \frac{(1+t)^{\dim L_{even} - a}}{(1-t)^b}
$$ 
and
$$
\pm \frac{1}{2^{|G_{odd}|}} \frac{(1+t)^{\dim L_{even} - c}}{(1-t)^d} \theta ,
$$
where $a$, $b$, $c$, $d$ are sums of dimensions of homogeneous components $L_\alpha$
for certain combinations of $\alpha\in G_{odd}$, 
with terms attaining the maximal possible value $\dim L_{odd}$ for $a$ and $c$, being
obtained by picking the summands 
$$\frac{1}{2} \frac{1}{(1 + t)^{\dim L_\alpha}}$$
and 
$$-\frac{1}{2} \frac{1}{(1 + t)^{\dim L_\alpha}} \theta$$
from the product of terms (\ref{super-symm}) for all $\alpha\in G_{odd}$.
The sum of all such terms is equal to
\begin{equation*}
\frac{(1+t)^{\dim L_{even} - \dim L_{odd}}}{2^{|G_{odd}|}} 
\sum_{i=0}^{|G_{odd}|} \binom{|G_{odd}|}{i} (-1)^i \theta^i = 
(1+t)^{\dim L_{even} - \dim L_{odd}} \Big( \frac{1 - \theta}{2} \Big)^{|G_{odd}|}.
\end{equation*}
From here part (ii) follows.

The general color case is yet more cumbersome, but similar. Combining 
(\ref{prod}), (\ref{color-skew}) and (\ref{color-symm}), we should care about
denominators in (\ref{color-symm}) vanishing for $t=-1$.
The corresponding summands in (\ref{color-symm}) obtained for $j=|\alpha|/2$ 
(recall that since $\alpha\in G_{odd}$, $|\alpha|$ is even) and are of the form 
$$
\frac{(-1)^i}{|\alpha|} \frac{1}{(1+t)^{\dim L_\alpha}} e^{i\alpha}, 
\quad i=0, \dots, |\alpha| - 1 .
$$
The maximal possible value of a power of $(1+t)$
in denominator is attained when we pick these summands from the product of terms 
(\ref{color-symm}) for all $\alpha\in G_{odd}$. The sum of all such terms is equal to
\begin{equation}\label{denom}
\frac{1}{(1+t)^{\dim L_{odd}}} \prod_{\alpha\in G_{odd}} \overline\chi(\alpha) .
\end{equation}
This could be compensated only by factors from (\ref{color-skew}) having
a power of $(1+t)$ in numerator, which are obtained for $j=0$.
The \textit{minimal} possible value $\dim L_0$ for a power of $(1+t)$ in numerator 
is attained when we pick the only possible summand $(1+t)^{\dim L_0}$ corresponding to 
$0 \in G$, and any summand \textit{not} containing a power of $(1+t)$ for all other 
$0 \ne \alpha\in G_{even}$ (that is, any summand with $j \ne 0$).
Therefore, the resulting expression is equal to the product of all terms (\ref{color-skew})
(or, what is equivalent, (\ref{skew-finite}))
for all $\alpha\in G_{even}$, minus product of terms with $j=0$ for at least one 
$\alpha\ne 0$:
\begin{equation}\label{numer}
(1+t)^{\dim L_0} \prod_{0\ne \alpha\in G_{even}} 
\Big( (1+te^\alpha)^{\dim L_\alpha} - 
      \frac{1}{|\alpha|} (1+t)^{\dim L_\alpha} \sum_{i=0}^{|\alpha| - 1} e^{i\alpha} \Big) .
\end{equation}
Comparison of the limiting behavior at $t \to -1$ of (\ref{denom}) and (\ref{numer})
completes the proof.
\end{proof}

\section{So, when it indeed vanishes?}\label{section-vanishes}

Finally, we arrive at our main result:

\begin{theorem}\label{theorem-main}
Let $L$ be a finite-dimensional color Lie algebra, $M$ be a finite-di\-men\-si\-o\-nal $L$-module,
and $B^*$ be the space of coboundaries in the corresponding color
Chevalley-Eilenberg complex.
\begin{enumerate}
\item
If $\chi(B^*)$ exists, then $\chi(L,M)$ exists and is equal to
$$
\begin{cases}
0,                               & L_{even} \ne 0  \\
\frac{\dim M}{2^{\dim L_{odd}}}, & L_{even} = 0.
\end{cases}
$$

\item 
If $\chi_{\mathbb Z_2} (B^*)$ exists and $\dim L_{even} \ge \dim L_{odd}$, then 
$\chi_{\mathbb Z_2} (L,M)$ exists and is equal to
$$
\begin{cases}
0,                                        & \dim L_{even} > \dim L_{odd}  \\
\big( \frac{1 - \theta}{2} \big)^{|G_{odd}|} \dim_{\mathbb Z_2} M, & \dim L_{even} = \dim L_{odd}.
\end{cases}
$$

\item
If $\chi_G(B^*)$ exists and $\dim L_0 \ge \dim L_{odd}$, then $\chi_G(L,M)$ exists 
and is equal to
$$
\begin{cases}
0,                                                        & \dim L_0 > \dim L_{odd}  \\
\prod_{0\ne \alpha\in G_{even}} (1 - e^\alpha)^{\dim L_\alpha}
\prod_{\alpha\in G_{odd}} \overline\chi(\alpha) \dim_G M, & \dim L_0 = \dim L_{odd}.
\end{cases}
$$
\end{enumerate}
\end{theorem}

\begin{proof}
Let us prove part (iii).
By Lemma \ref{complex}(iii), $\chi_G(\bigwedge\nolimits_\varepsilon^* L)$ exists.
By Lemma \ref{lemma-dimension}(iii), $\chi_G(C^*(L,K))$ exists and coincides with
$\chi_G(\bigwedge\nolimits_\varepsilon^* L)$.
Then by Lemma \ref{lemma-dimension}(ii), $\chi_G(C^*(L,M))$ exists and is equal to
$\chi_G(C^*(L,K)) \dim_G M$.
Finally, by Lemma \ref{principle}, $\chi_G(L,M)$ exists and coincides with 
$\chi_G(C^*(L,M))$, what finishes the proof.

Parts (i) and (ii) are proved by the same implications, replacing
$\chi_G$ by $\chi$ and by $\chi_{\mathbb Z_2}$, respectively.
\end{proof}

\begin{corollary}
Let $L$ be a finite-dimensional color Lie algebra, $M$ be a finite-di\-men\-si\-o\-nal
$L$-module such that cohomology $H^*(L,M)$ is finite-dimensional.
\begin{enumerate}
\item If $L_{even} \ne 0$, then $\chi(L,M) = 0$.
\item If $\dim L_{even} > \dim L_{odd}$, then $\chi_{\mathbb Z_2}(L,M) = 0$.
\item If $\dim L_0 > \dim L_{odd}$, then $\chi_G(L,M) = 0$.
\end{enumerate}
\end{corollary}

\begin{proof}
Let us, again, prove part (iii).
$\chi_G(L,M)$, being a limit of a finite sum, obviously exists, and in view of
Theorem \ref{theorem-main},
it is sufficient to prove the existence of $\chi_G(B^*)$.
As ignoring the fixed number of initial terms of a series does not affect its convergence,
for this purpose we may ignore all initial terms up to the degree of the highest 
nonzero cohomology,
and assume $H^* = 0$ and $Z^* = B^*$. Then by Lemma \ref{lemma-one-plus-t},
$\chi_G(C^*,t) = (1+t) \chi_G(B^*,t)$, and by Lemma \ref{lemma-dimension}(iii), 
$\chi_G(C^*,t) = \chi_G(\bigwedge\nolimits_\varepsilon^* L, t)$.
Hence, for the Poincar\'e series of $B^*$ the same counting of powers of $(1+t)$ as for 
$\bigwedge\nolimits_\varepsilon^* L$ in the proof of Lemma \ref{complex} applies,
but shifted by $-1$, and the necessary and sufficient condition for the existence of 
$\chi_G(B^*)$ is $\dim L_0 \ge \dim L_{odd} + 1$, what finishes the proof.

And, again, parts (i) and (ii) are absolutely similar.
\end{proof}

\section{Another method to the rescue?}\label{section-borel}

So, after all, the Abel summation turned out to be not so powerful for our purposes:
according to Lemma \ref{complex}, the Euler-Poincar\'e characteristic of the (co)chain
Chevalley-Eilenberg complex does not exist in many cases.
Perhaps, we should try another method?

For the basic facts concerning summation of divergent series
we refer to the classical book \cite{hardy}. However, since that time some modern
terminology has been adopted, which we briefly remind here.

A summation method $(S)$, assigning to some (divergent) series with real coefficients 
$\sum_{n=0}^\infty a_n$ a real value, denoted as $(S) \sum_{n=0}^\infty a_n$, is called 
\textit{regular}, if any convergent series is summable by this method to its limit, 
\textit{additive} if from summability of series $\sum_{n=0}^\infty a_n$ and 
$\sum_{n=0}^\infty b_n$ follows summability of their linear combination 
$\sum_{n=0}^\infty (\lambda a_n + \mu b_n)$ for any $\lambda, \mu \in \mathbb R$, and
$$
(S) \sum_{n=0}^\infty (\lambda a_n + \mu b_n) = 
\lambda \cdot (S) \sum_{n=0}^\infty a_n + \mu \cdot (S) \sum_{n=0}^\infty b_n
$$
holds,
\textit{multiplicative} if from summability of series $\sum_{n=0}^\infty a_n$ and 
$\sum_{n=0}^\infty b_n$
follows summability of their (Cauchy) product 
$\sum_{n=0}^\infty c_n = \sum_{n=0}^\infty \big (\sum_{i=0}^n a_i b_{n-i} \big)$, and
$$
(S) \sum_{n=0}^\infty c_n = (S) \sum_{n=0}^\infty a_n \cdot (S) \sum_{n=0}^\infty b_n
$$
holds,
and \textit{left-translative} if from summability of series $\sum_{n=0}^\infty a_n$ follows
summability of the series $\sum_{n=1}^\infty a_n$, and
$$
(S) \sum_{n=1}^\infty a_n = -a_0 + (S) \sum_{n=0}^\infty a_n
$$
holds.

Let $(S)$ be a summation method.
If we can sum, $G$-componentwise, the series (\ref{sum}) by this method,
we will call the resulting sum (which is still an element of $\mathbb R[G]$)
an \textit{Euler-Poincar\'e characteristic in the sense of (S)} and denote it by $\chi_G^{(S)}$
(so, the Euler-Poincar\'e characteristic defined in \S \ref{section-char} becomes 
the Euler-Poincar\'e characteristic in the sense of Abel).
To be consistent with the usual (finite) notion of the Euler-Poincar\'e characteristic,
the method should be regular.

In this general situation, the nice conceptual link to Poincar\'e series is lost.
However, essentially the same reasonings as in 
\S\S \ref{section-char} and \ref{section-vanishes} could be used to prove similar results.

\begin{lemma}
Let $(S)$ be a linear and left-translative summation method,
$C^*$ be a color cochain complex all whose graded components are finite-dimensional,
and $B^*$ be the space of its coboundaries.
If both $\chi_G^{(S)}(C^*)$ and $\chi_G^{(S)}(B^*)$ exist,
then $\chi_G^{(S)}(H(C^*))$ exists and coincides with $\chi_G^{(S)}(C^*)$.
\end{lemma} 

\begin{proof}
By linearity, (\ref{exact-1}) and (\ref{exact-2}) imply respectively:
\begin{align*}
\chi_G^{(S)}(C^*) &= \chi_G^{(S)}(Z^*) + \chi_G^{(S)}(B^*[1])  \\
\chi_G^{(S)}(Z^*) &= \chi_G^{(S)}(B^*) + \chi_G^{(S)}(H(C^*)) .
\end{align*}
By left-translativity, 
$$
\chi_G^{(S)}(B^*[1]) = \dim_G B_0 - \chi_G^{(S)}(B^*) = -\chi_G^{(S)}(B^*) ,
$$
so the first equality implies existence of $\chi_G^{(S)}(Z^*)$, then the second one
implies existence of $\chi_G^{(S)}(H(C^*))$, and together they imply 
$\chi_G^{(S)}(H(C^*)) = \chi_G^{(S)}(C^*)$.
\end{proof}

\begin{lemma}\label{lemma-summab}
Let $(S)$ be a regular and multiplicative summation method, and $L$ be a finite-dimensional 
color Lie algebra. If for any $\alpha\in G_{odd}$ and any $0 \le n < |\alpha|$, 
the series $$\sum_{i=0}^\infty \binom{\dim L_\alpha + |\alpha|i + n - 1}{|\alpha|i + n}$$
is (S)-summable, then $\chi_G^{(S)} (\bigwedge\nolimits_\varepsilon^* L)$ exists.
If additionally, $L_0 \ne 0$, then $\chi_G^{(S)} (\bigwedge\nolimits_\varepsilon^* L) = 0$.
\end{lemma}

\begin{proof}
By multiplicativity, 
\begin{equation}\label{twist-borel}
\chi_G^{(S)} (V \otimes W) = \chi_G^{(S)}(V) \chi_G^{(S)}(W)
\end{equation}
for any two $G$-color vector spaces, with existence of the Euler-Poincar\'e characteristics
on the right-hand side implies the existence of Euler-Poincar\'e characteristics on the 
left-hand side.
Note that, unlike (\ref{twist}), this is no longer a sole consequence of multiplication rule 
in a group ring, but follows from the rule of series multiplication and postulated 
multiplicativity of $(S)$.
Twisting of tensor product (like in \S \ref{section-color}) obviously does not
affect the left-hand side of (\ref{twist-borel}), so the same equality holds for the twisted 
tensor product $\widehat\otimes$.

Now, (\ref{decomp}) implies, in a fashion similar to (\ref{prod}):
$$
\chi_G^{(S)} (\bigwedge\nolimits_\varepsilon^* L) =
\prod_{\alpha\in G_{even}} \chi_G^{(S)} (\bigwedge\nolimits^*(L_\alpha))
\prod_{\alpha\in G_{odd}}  \chi_G^{(S)} (\bigS^*(L_\alpha)) .
$$
By regularity, and similarly to (\ref{skew-finite}):
\begin{equation*}
\chi_G^{(S)} (\bigwedge\nolimits^*(L_\alpha)) = 
\sum_{n=0}^{\dim L_\alpha} (-1)^n \dim \bigwedge\nolimits^n(L_\alpha) \> e^{n\alpha} =
\sum_{n=0}^{\dim L_\alpha} (-1)^n \binom{\dim L_\alpha}{n} \> e^{n\alpha} = 
(1 - e^\alpha)^{\dim L_\alpha}
\end{equation*}
what vanishes for $\alpha = 0$, so to establish the desired
conclusions it is sufficient to prove the existence of 
$\chi_G^{(S)} (\bigS^*(L_\alpha))$'s.

Similarly to (\ref{symm}), we have:
\begin{equation*}
\chi_G^{(S)}(\bigS^*(L_\alpha)) = 
(S) \sum_{n=0}^\infty (-1)^n \dim \bigS^n (L_\alpha) \> e^{n\alpha} =
(S) \sum_{n=0}^\infty (-1)^n \binom{\dim L_\alpha + n-1}{n} \> e^{n\alpha} ,
\end{equation*}
what implies, similarly to (\ref{color-symm}) 
(but we cannot use multisections anymore, as denominator
in the corresponding expression vanishes):
$$
\chi_G^{(S)} (\bigS^*(L_\alpha)) = 
\sum_{n=0}^{|\alpha| - 1} (-1)^n \Big( 
(S) \sum_{i=0}^\infty \binom{\dim L_\alpha + |\alpha|i + n - 1}{|\alpha|i + n}
\Big) e^{n\alpha} ,
$$
but these sums exist by the assumption.
\end{proof}

The just proved Lemma is very far from being the strongest result in this direction:
we were lazy and treated each $\chi_G^{(S)} (\bigS^*(L_\alpha))$ separately, while, 
intermixing the powers $e^{n\alpha}$ for different $\alpha\in G_{odd}$ could lead to 
more relaxed sufficient conditions.
Our purpose was merely to demonstrate that, choosing a strong enough summation method,
we may achieve existence (and vanishing in the generic case) of the 
Euler-Poincar\'e characteristic in the sense of this method. 
Note that comparing Lemmata \ref{complex} and
\ref{lemma-summab} reveals that any such method would be necessarily
incompatible with the Abel summation, as there are situations in which the 
Euler-Poincar\'e characteristic does not vanish in one sense and vanishes in another.

\begin{theorem}\label{theorem-borel}
Let $L$ be a finite-dimensional color Lie algebra, $M$ be a finite-di\-men\-si\-o\-nal $L$-module,
$B^*$ be the space of coboundaries in the corresponding color Chevalley-Eilenberg complex,
and $(S)$ be a regular, additive, multiplicative and left-translative summation method.
Suppose that $\chi_G^{(S)}(B^*)$ exists and that the summability assumption of 
Lemma \ref{lemma-summab} holds. Then $\chi_G^{(S)}(L,M)$ exists. If, additionally, 
$L_0 \ne 0$, then $\chi_G^{(S)}(L,M) = 0$.
\end{theorem}

\begin{proof}
The proof repeats the proof of Theorem \ref{theorem-main}, with the difference that, again, instead
of appealing to Lemma \ref{lemma-dimension}, we use equality (\ref{twist-borel})
and an obvious equality $\chi_G^{(S)}(M) = \dim_G M$ which follows from the regularity.
\end{proof}

Note, however, that we cannot derive an analog of Corollary to Theorem 1.

\section{Further questions}\label{section-questions}

\subsection*{1}
Could we get rid of the condition of existence of $\chi_G(B^*)$ in Theorem \ref{theorem-main}
or Theorem \ref{theorem-borel}?
(From the proof of Lemma \ref{lemma-one-plus-t} it is clear that this is equivalent to 
the existence of $\chi_G(Z^*)$).
This is so, for example, for all color Lie algebras of ordinary dimension $3$,
as follows from the results of  
\cite{piont-s}\footnote[2]{There is a misprint in \cite{piont-s}: in the Table 1,
for algebra 1, in the case $\mu = -1$, the (ordinary) Poincar\'e series should be
$1 + z + z^2 + z^3$. This is rectified in the arXiv version of the paper.}.

\subsection*{2}
What happens when we are unable to apply the 
passing-to-the-underlying-cochain-complex argument? This applies both to the 
finite-dimensional case where the Euler-Poincar\'e characteristic of the underlying 
cochain complex does not exist (i.e., when conditions on dimensions in Theorem \ref{theorem-main} and its Corollary are not satisfied), and to the infinite-dimensional case
(i.e., when either a Lie superalgebra, or its module, or both, are infinite-dimensional).
In such situations, it seems that we inevitably have to avoid tricks and start to do a real job
by looking at cohomology itself.
Perhaps, a not very deep conjecture would be:
\begin{conjecture}
Let $L$ be a (possibly infinite-dimensional) $G$-color Lie algebra with a ``large enough''
zero part, and $M$ be a (possibly infinite-dimensional) $L$-module 
such that cohomology $H^*(L,M)$ is finite-dimensional in each degree. 
Then $\chi_G(L,M)$ exists and vanishes.
\end{conjecture}

The question seems to be nontrivial (and apparently never investigated)
already in the very particular cases: for an infinite-dimensional ordinary Lie algebra 
$L$, and when the whole cohomology $H^*(L,M)$ of a color Lie algebra $L$ is 
finite-dimensional (i.e., vanishes for all but the finite number of degrees;
the Euler-Poincar\'e characteristic then clearly exists).

There are numerous examples of infinite-dimensional Lie algebras of vector fields and 
close to them, satisfying  the condition of the conjecture 
(see the classical book \cite{fuchs} and a more recent compendium of
results in \cite[Chapter 3]{ff}), and it seems that 
in all the cases where cohomology is known, the conjecture holds.

Still, may be this conjecture is too optimistic: it implies, in particular, that the radius
of convergence of the Poincar\'e series of any algebra which satisfies the condition
of the conjecture, is equal to $1$. So, it would be more cautious to ask first
whether this is indeed so. However, if we no longer deal with the Euler-Poincar\'e characteristic,
but with entire Poincar\'e series, there is no need to assume restriction on dimensions
of the parts of an algebra.
Thus, a more appropriate formulation of this cautious question, 
which we learned from Dmitri Piontkovski, would be:
does there exist a $G$-color Lie algebra $L$ and an $L$-module $M$ such that
cohomology $H^*(L,M)$ is finite-dimensional in each degree and the radius of convergence 
of some $G$-component of $\chi_G(L,M,t)$ is $<1$? (The latter condition may be 
reformulated as follows: some $G$-component of $\dim_G H^n(L,M)$ grows exponentially with 
$n$).

Another, significantly weakened, form of this conjecture may be 
obtained by observing that in all the established cases the Euler-Poincar\'e characteristic, 
even when it does not vanish, involves elements $\overline\chi(\alpha)$ which are zero divisors
in $\mathbb Z[G]$. This leads to the following

\begin{conjecture}
Let $G$ be an abelian group with torsion, $L$ be a (possibly infinite-dimensional) $G$-color Lie 
algebra with a ``large enough'' zero part, and $M$ be a (possibly infinite-dimensional) $L$-module 
such that cohomology $H^*(L,M)$ is finite-dimensional in each degree. 
Then $\chi_G(L,M)$ exists and is a zero divisor in $\mathcal C[G]$.
\end{conjecture}

\subsection*{3}
Could reasonings based on the Hochschild-Serre spectral sequence, which are redundant
in the ordinary Lie algebra case, be useful in the super or color case?

\subsection*{4}
Another approach to the Euler-Poincar\'e characteristic may be pursued via 
Poincar\'e duality. It seems to be interesting to investigate
what would be a Poincar\'e duality analog for a Lie superalgebra cohomology (if any)?

\subsection*{5}
In the true Eulerian spirit, instead of using Lemma \ref{multisec}, one may wish to 
write $\chi_G(W, t) = \chi(W, te^\alpha)$. In particular, (\ref{symm}) could be
rewritten in this way as
$$
\chi_G(\bigS^*(L_\alpha), t) = \frac{1}{(1 - te^\alpha)^{\dim L_\alpha}},
$$
and then the whole Poincar\'e series of the color Chevalley-Eilenberg chain complex
acquires the attractive shape
\begin{equation}\label{spec}
\frac
{\prod_{\alpha \in G_{even}} (1 + te^\alpha)^{\dim L_\alpha}}
{\prod_{\alpha \in G_{odd}} (1 - te^\alpha)^{\dim L_\alpha}}
\end{equation}
which, perhaps, gives a chance to avoid more cumbersome reasonings based on multisection of 
series. To make any sense of this, one should consider Poincar\'e series and 
Euler-Poincar\'e characteristic with values in something like the ring of quotients of 
the completion of $\mathbb R[G]$, and allow limiting processes which are compatible with 
the Euler-Poincar\'e characteristic as understood in this note.

Compare (\ref{spec}) with the formula at the bottom of p. 604 in \cite{kang}, which
expresses the dimension of the $\varepsilon$-symmetric power of a color vector space
(that is, odd and even parts being swapped).
Note that the situation here, however, is somewhat different: the group $G$ has a torsion, and hence, 
as was already noted,  elements $\overline\chi(\alpha)$ are zero divisors 
in $\mathbb Z[G]$, so the straightforward notion of the ring of quotients is not applicable.


\section*{Acknowledgements}

Thanks are due to Rutwig Campoamor-Stursberg for triggering my interest in the subject 
and pointing out some references, to Dmitri Piontkovski for an interesting discussion and 
informing about misprints in \cite{piont-s}, 
to Max Alekseyev for advice to use multisections in Lemma \ref{multisec},
and to Dimitry Leites for interest and useful remarks.
On the final stage of prepartion of this manuscript, I was supported by the grant ERMOS7 
co-funded by the Estonian Science Foundation and Marie Curie Actions.


\begin{thebibliography}{Bou}

\bibitem[BL]{berger-l} C. Berger and T. Leinster,
\emph{The Euler characteristic of a category as the sum of a divergent series},
Homology, Homotopy and Appl. \textbf{10} (2008), 41--51; \textsf{arXiv:0707.0835}.

\bibitem[Bor]{borcherds} R.E. Borcherds, \emph{Sporadic groups and string theory},
First European Congress of Mathematics, Vol. I: Invited Lectures, Birkh\"auser, 1994, 
411--421.

\bibitem[Bou]{bourbaki} N.~Bourbaki, 
\emph{Alg\`ebre. Chapitre X. Alg\`ebre homologique}, Masson, Paris, 1980.

\bibitem[CE]{ce} C. Chevalley and S. Eilenberg, 
\emph{Cohomology theory of Lie groups and Lie algebras}, 
Trans. Amer. Math. Soc. \textbf{63} (1948), 85--124.

\bibitem[FF]{ff} B.L. Feigin and D.B. Fuchs, 
\emph{Cohomologies of Lie groups and Lie algebras}, Encyclopaedia of Mathematical Sciences 
Vol. 21, Lie Groups and Lie Algebras II, Springer, 2000, 127--215.

\bibitem[F]{fuchs} D.B. Fuchs, \emph{Cohomology of Infinite-Dimensional Lie Algebras},
Consultants Bureau, New York, 1986.


\bibitem[Go]{goldberg} S.I. Goldberg, \emph{On the Euler characteristic of a Lie algebra},
Amer. Math. Monthly \textbf{62} (1955), 239--240.

\bibitem[Gr]{gruson} C. Gruson, 
\emph{Finitude de l'homologie de certains modules de dimension finie sur 
une super alg\`ebre de Lie}, Ann. Inst. Fourier \textbf{47} (1997), 531--553.

\bibitem[H]{hardy} G.H. Hardy, \emph{Divergent Series}, Oxford Univ. Press, 1949.

\bibitem[Ka]{kang} S.-J. Kang, 
\emph{Graded Lie superalgebras and the superdimension formula}, 
J. Algebra \textbf{204} (1998), 597--655.

\bibitem[Ko]{koszul} J.-L. Koszul, \emph{Homologie et cohomologie des alg\`ebres de Lie},
Bull. Soc. Math. France \textbf{78} (1950), 65--127; reprinted in 
\emph{Selected Papers}, World Scientific, 1994, 7--69.

\bibitem[M1]{malliavin-nilp} M.-P. Malliavin-Brameret, 
\emph{Cohomologie d'alg\`ebres de Lie nilpotentes, et caract\'eristiques 
d'Euler-Poincar\'e}, Bull. Sci. Math. \textbf{100} (1976), 269--287.

\bibitem[M2]{malliavin-solv} \bysame, 
\emph{Caract\'erisitiques d'Euler-Poincar\'e d'alg\`ebres de Lie r\'esolubles},
Comptes Rendus Acad. Sci. Paris \textbf{284} (1977), 1487--1488.

\bibitem[PS]{piont-s} D. Piontkovski and S. Silvestrov, 
\emph{Cohomology of $3$-dimensional color Lie algebras}, 
J. Algebra \textbf{316} (2007), 499--513; \textsf{arXiv:math/0508573}.

\bibitem[P]{pirashvili} T. Pirashvili, 
\emph{The Euler-Poincar\'e characteristic of a Lie algebra},
J. Lie Theory \textbf{8} (1998), 429--431.

\bibitem[R]{riordan} J. Riordan, \emph{Combinatorial Identities}, Wiley, 1968.

\bibitem[S]{scheunert} M. Scheunert, \emph{Generalized Lie algebras}, 
J. Math. Phys. \textbf{20} (1979), 712--720. 

\bibitem[SZ]{scheunert-zhang} \bysame{ }and R.B. Zhang, 
\emph{Cohomology of Lie superalgebras and of their generalizations},
J. Math. Phys. \textbf{39} (1998), 5024--5061; \textsf{arXiv:q-alg/9701037}.

\end{thebibliography}
\end{document}